\documentclass[11pt]{article}
\usepackage{amsmath,amssymb,amsfonts,amsthm}
\usepackage{mathrsfs}
\usepackage{indentfirst}
\usepackage[pdftex]{color,graphicx}
\usepackage{subfig}
\usepackage[pdftex,bookmarks,unicode,colorlinks]{hyperref}
\usepackage{enumerate} 
\usepackage[numbers]{natbib} 
\usepackage{yfonts} 
\usepackage{caption} 
\usepackage{bm}

\definecolor{awesome}{rgb}{1.0, 0.13, 0.32}
\definecolor{bronze}{rgb}{0.8, 0.5, 0.2}
\definecolor{burntorange}{rgb}{0.8, 0.33, 0.0}
\definecolor{blue(ncs)}{rgb}{0.0, 0.53, 0.74}
\theoremstyle{plain}
\newtheorem{theorem}{Theorem}[section]

\theoremstyle{remark}
\newtheorem{remark}[theorem]{Remark}
\newtheorem{example}[theorem]{Example}

\theoremstyle{definition}

\numberwithin{equation}{section} 

\begin{document}

\title{\bf An Optimal Control Method to Compute the Most Likely Transition Path for Stochastic Dynamical Systems with Jumps }

 \author{\bf\normalsize{
Wei Wei$^{1,}$\footnotemark[2],
Ting Gao$^{1,}$\footnotemark[3],
Xiaoli Chen$^{2,}$\footnotemark[1],
and Jinqiao Duan$^{3,}$\footnotemark[4]
}\\[10pt]
\footnotesize{$^1$Center for Mathematical Sciences, Huazhong University of Science and Technology,} \\
\footnotesize{Wuhan, Hubei 430074, China.} \\[5pt]
\footnotesize{$^2$ Department of Mathematics, National University of Singapore, 119077, Singapore.}\\
\footnotesize{$^3$ Department of Applied Mathematics \& Department of Physics, Illinois Institute of Technology, Chicago, IL 60616, USA.}
}

\footnotetext[2]{Email: \texttt{weiw16@hust.edu.cn}}
\footnotetext[3]{Email: \texttt{tgao0716@hust.edu.cn}}
\footnotetext[1]{Email: \texttt{xlchen@nus.edu.sg}}
\footnotetext[4]{Email: \texttt{duan@iit.edu}}
\footnotetext[1]{is the corresponding author}

\date{March 28, 2022}
\maketitle
\vspace{-0.3in}

\begin{abstract}
Many complex real world phenomena exhibit abrupt, intermittent or jumping behaviors, which are more suitable to be described by stochastic differential equations under non-Gaussian L\'evy noise. Among these complex phenomena, the most likely transition paths between metastable states are important since these rare events may have a high impact in certain scenarios. Based on the large deviation principle, the most likely transition path could be treated as the minimizer of the rate function upon paths that connect two points. One of the challenges to calculate the most likely transition path for stochastic dynamical systems under non-Gaussian L\'evy noise is that the associated rate function can not be explicitly expressed by paths. For this reason, we formulate an optimal control problem to obtain the optimal state as the most likely transition path. We then develop a neural network method to solve this issue. Several experiments are investigated for both Gaussian and non-Gaussian cases.  

\end{abstract}


{\bf Lead paragraph:}

\textbf{ Random fluctuations in dynamical systems may cause unexpected rare events which can be characterized through transitions between metastable states. In this paper, we investigate the transition behaviors of a stochastic dynamical system under non-Gaussian L\'evy noise, which is a more general noise than Gaussian noise that is widely used in models of climate, thermal dynamics and engineering. First, the rate function from the large deviation principle is applied to estimate the probability of transitions through a small tube of paths, which connect two given points. A constrained minimization problem is then set up for obtaining the minimizer of the rate function upon paths connecting two given points, which is the most likely transition path. Due to the implicit form of the rate function with respect to such paths, the existing numerical methods for Gaussian dynamics, for example, the minimum action method, can not be directly applied. As a result, 
we reformulate the constrained minimization problem as an optimal control problem and the optimal state is the desired most
likely transition path. However, the existing methods in the optimal control problem can not be directly applied to our case. We then construct two neural networks to approximate the state and control with a loss function which makes the state and control satisfy the constraints, and further, the cost functional attains its minimum value. Training the neural networks with such loss function, we obtain the desired most likely transition path. Numerical experiments on the Maier Stein system under Gaussian and non-Gaussian L\'evy noise show the validation of our method.}


\section{ Introduction}

Random fluctuations in dynamical systems may lead to rare events with high impact under certain circumstances\cite{duan_book,Lucarini, Schweiger}. One type of these rare events, transition pathways between metastable states, are important research topics in stochastic dynamical systems\cite{Kuehn, Raser,xu2020}. However, this kind of noise induced transition phenomena will not possibly happen in deterministic systems. 
In recent years, L\'evy noise has been playing an increasingly important role in models of climate\cite{ Ganopolski2002}, biology\cite{xu2013,zheng2016,chen2019}, engineering\cite{Sanchez2015} and financial markets\cite{yuan2022}.
One of the L\'evy noise induced dynamical systems, the L\'evy flight, plays an important role in anomalous diffusion\cite{PR195, PR339}, optimized search\cite{PNAS111, Cambridge2011} and stability analysis\cite{EPL72, PRE75, PRE67, PRE59, JSM2010}.
Such a growing amount of research shows that non-Gaussian L\'evy process is more appropriate to be utilized as the intrinsic noise in various complex systems. As a result, this inspires us to investigate the transition pathways of stochastic dynamical systems under non-Gaussian L\'evy noise.

The large deviation principle is a powerful tool in tackling transition behaviors. In general, the large deviation principle measures how small the probability of rare events is through the rate function. It also generalized the concept of potential of the gradient system to quasi-potential of the non-gradient system, which is the minimum value of the rate function defined on paths connecting two points. The associated minimizer is regarded as the most likely transition path in a wide range of literature \cite {heymann2008, zhou2008, weinan2004,hu2021}. However, as far as we know, the major part of literature studying transition pathways is for stochastic dynamical systems under Gaussian noise, that is, the stochastic differential equations driven by Brownian motion. For the dynamics of stochastic dynamical systems under non-Gaussian L\'evy noise, the literature is relatively limited.

Several methods to calculate the most likely transition path for stochastic dynamical systems under Gaussian noise involve two aspects: solving an Euler-Lagrange equation or a constrained optimal control problem. 
As for the former one, the minimum action method(MAM) is a widely used numerical method in finding the most likely transition path based on the large deviation principle. A constrained minimization problem was proposed to compute the most likely transition path for stochastic differential equations with multiple metastable states~\cite{weinan2004}. The basic idea of MAM is to find the minimizer of the associated rate function and take this minimizer as the desired most likely transition path. 
 The geometric minimum action method(gMAM)\cite{heymann2008} improved the MAM through reparametrization and provided a series of theoretical results on the most likely transition path. The adaptive minimum action method(aMAM) \cite{zhou2008,sun2018} further improved the efficiency of MAM through a moving mesh method. The finite element method was employed in \cite{wan2018} and widened the scope of applications of the aMAM.
However, the MAM usually tackles the problem that the associated rate function has an explicit form expressed by paths. When it comes to the case involving a L\'evy noise, as the associated rate function does not have such kind of explicit form, the MAM is not applicable in general.

While for the latter aspect, the constrained minimization problem associated with computation of the most likely transition path can be naturally regarded as an optimal control problem. There is a huge amount of literature in the field of optimal control theory, involving time or space-marching \cite{WKang2021, WKang2015}, neural network \cite{WKang2021,Li2017,Li2018}, stochastic process\cite{HanJ2018}, Hopf formula\cite{Darbon2016} and et al. For the stochastic dynamical systems under Gaussian noise, the most likely transition path is directly the solution of a common optimal control problem as \cite[Equation (2.1)]{WKang2021} and similar work may be found in \cite{Todorov2011, YGao2020}. However, for stochastic dynamical systems under non-Gaussian L\'evy noise, the control involves space variable as a free variable which makes itself not a direct application of existing methods in the field of optimal control problems.

In recent years, the machine learning method has greatly influenced a lot of areas such as numerical computing. It sheds light on solving various kinds of partial differential equations effectively \cite{HanJ2018,Maziar2019,chen2021}. This gives us a new viewpoint to compute the minimizer of an optimal control problem through constructing neural networks. Therefore we leverage the deep learning method to solve our problem on the most likely transition path.

In this paper, we consider the transition problem of the stochastic dynamical systems under L\'evy noise with exponentially light jumps, which is described by stochastic differential equations. Our aim is to find the most likely transition path between metastable states of the stochastic differential equations as the noise intensity tends to $0$. Based on the large deviation principle, we formulate the problem as an optimal control problem and the optimal state is the desired most likely transition path. As the existing methods in the optimal control problem can not be directly applied to our case, we then develop a neural network method to approximate the state and control. 
 Similar as the idea of the soft constraints method\cite{lu2021physics}, we construct a loss function which makes the state and control satisfy the constraints, and further, we obtain the desired most likely transition path when the cost functional attains its minimum value. 

As a validation for our method, we first apply it to the Maier-Stein system under Gaussian noise with two different parameters, and compare our numerical results with the most likely transition path that is computed  through the adaptive minimum action method. Correspondingly, we also compute the results for the Maier-Stein system under non-Gaussian L\'evy noise with exponentially light jumps, and  uncover some intrinsic relationship between two different types of noises. To our knowledge, this is the first numerical result on the most likely transition path to stochastic differential equations driven by non-Gaussian L\'evy noise in the context of the large deviation principle.

This paper is arranged as follows. In section \ref{B}, we briefly introduce the large deviation principle and the constrained minimization problem in finding the most likely transition path. In section \ref{ML} we formulate the constrained minimization problem as an optimal control problem. We then construct neural networks and a proper loss function to solve the optimal control problem. Finally, numerical results on the Maier Stein system under Gaussian noise and non-Gaussian L\'evy noise with exponentially light jumps are presented in section \ref{NE}.

\section{Background}\label{B}

We consider the following $d$-dimensional stochastic differential equation under non-Gaussian L\'evy noise:
\begin{equation}\label{maineq}
     d X_t^{\varepsilon}=b\left(X_t^{\varepsilon}\right) d t+\varepsilon d \tilde{L}_t^{\varepsilon},\quad X_0^\epsilon=x_1.
\end{equation}
Function $b: \mathbb{R}^d \to \mathbb{R}^d$ is the drift term.
For every $\varepsilon>0$, the non-Gaussian L\'evy process $(\tilde{L}_{t}^{\varepsilon})_{t \geq 0}$ is given by
\begin{equation}\label{Lt}
\tilde{L}_{t}^{\varepsilon}=\int_{0}^{t} \int_{\mathbb{R}^d\backslash \{0\}} z \tilde{N}^{\frac{1}{\varepsilon}}(d s, d z),
\end{equation}
where $\tilde{N}^{\frac{1}{\varepsilon}}$ is a compensated Poisson random measure defined on a given complete probability space $(\Omega, \mathcal{F}, \mathbb{P})$ with compensator $\varepsilon^{-1} d s \otimes \nu$. The measure $\nu$ has the form
\begin{equation}
\nu(d z)=e^{-|z|^{\gamma}} d z.
\end{equation}
 Here, the parameter $\gamma>1$. According to the L\'evy-Ito's decomposition, every L\'evy process can be expressed as the sum of a continuous part and a part expressible as a compensated sum of independent jumps\cite{sato}. For the L\'evy process $\tilde{L}_t^\epsilon$, it is a pure jump L\'evy process and thus expressed as the integral of the compensated random Poisson measure $\tilde{N}^{\frac{1}{\epsilon}}$.

The characteristic function for the process $\tilde{L}_{t}^{\varepsilon}$ is
\begin{equation}
    \mu_t(x)=\exp \left(\frac{t}{\epsilon}\int_{\mathbb{R}^d\backslash \{0\}} (e^{i(\xi,x)-|\xi|^\gamma}-e^{-|\xi|^{\gamma}}) d\xi\right).
\end{equation}

Comparing to the Langevin approach \cite{PR195, PR339, PRL73, PRE50}, the drift term $b$ represents the external force field. The L\'evy process $\tilde{L}_{t}^{\varepsilon}$ can be regarded as a continuous time random walk with waiting time obeying an exponential distribution with mean $t \nu(\mathbb{R}^d\backslash \{0\})$ and the step size obeying the distribution $\nu / \nu(\mathbb{R}^d\backslash \{0\})$ up to time $t$. 

According to the characteristic function of L\'evy process $\tilde{L}_{t}^{\varepsilon}$, the mean square displacement for the driven noise $\epsilon \tilde{L}_{t}^{\varepsilon}$ is 
\begin{equation}\label{MSD}
    \mathbb{E}\big(\epsilon \tilde{L}_{t}^{\varepsilon}\big)^2=\epsilon t \int_{\mathbb{R}^d\backslash \{0\}} |z|^2 e^{-|z|^{\gamma}} d z= K\epsilon t, 
\end{equation}
where the constant $K$ is determined by $\gamma$.

We want to investigate the transition behavior of equation \eqref{maineq} when the noise intensity $\epsilon$ is fairly small. The  most likely transition path is then introduced based on the large deviation principle in the following subsections.

\subsection{The large deviation principle}
For every $T>0$ and $x_1 \in \mathbb{R}^{d}$, let $\mathbb{D}(0,T)$ denote the Skorokhod space on $[0,T]$ equipped with the Skorokhod metric $\rho_T$, which is the space of real-valued right continuous functions with left limits. Let $C([0,T])$ denote the continuous functions on $[0,T]$. Define the deterministic mapping $F$ by $\varphi=F(g)$, where $\varphi$ is the solution of 
\begin{equation}\label{eq_phi}
    \varphi(t)=x+\int_{0}^{t} b(\varphi(s)) d s+\int_{0}^{t} \int_{\mathbb{R}^d\backslash \{0\}} z(g(s, z)-1) \nu(d z) ds.
\end{equation}
In \cite[Theorem 3.7]{Bud2013}, it was proved that equation \eqref{eq_phi} has a unique solution for every positive measurable function $g$ satisfying 
\begin{equation}\label{condg}
    \int_{0}^{T} \int_{\mathbb{R}^{d}\backslash \{0\}}(g(s, z) \ln g(s, z)-g(s, z)+1) \nu(d z) d s < \infty.
\end{equation}
We know that $X^{\varepsilon} $ obeys the large deviation principle \cite{de2021kramers} in the Skorokhod space $\mathbb{D}(0, T)$ with a good rate function $S_T: C[0,T] \longrightarrow$ $[0, \infty]$, given by
\begin{equation}\label{rf}
    S_T(\varphi) \triangleq  \inf_{\varphi=F(g)} \left\{\int_{0}^{T} \int_{\mathbb{R}^{d}\backslash \{0\}}(g(s, z) \ln g(s, z)-g(s, z)+1) \nu(d z) d s \right\},
\end{equation}
with the convention that $\inf \emptyset =\infty$. Here function $g$ is a positive measurable function on $[0, T]\times \mathbb{R}^d \backslash \{0\}$. Specifically, the distribution of solution $X^\epsilon$ satisfies
\begin{equation}\label{VaraLDP}
   -\inf_{\varphi \in A^\circ} S_T(\varphi) \leq \liminf_{\epsilon \to 0} \epsilon \log \mathbb{P}(X^\epsilon \in A) \leq \limsup_{\epsilon \to 0} \epsilon \log \mathbb{P}(X^\epsilon \in A) \leq -\sup_{\varphi \in \bar{A}} S_T(\varphi).
\end{equation}

As $S_T$ is a good rate function, the level set $\{S_T \leq s\} $ is compact in $\mathbb{D}(0, T)$ for every $s \geq 0$.  Then, according to \cite[Theorem 3.3]{FW12}, the Varadhan’s description of the large deviation principle \eqref{VaraLDP} is equivalent to the Freidlin-Wentzell's estimates.
 For small $\delta$, the Freidlin-Wentzell's estimates enable us to quantify the probability that $X^\epsilon$ stays in a $\delta$-tube of a path $\psi$ though the following asymptotic relation \cite{weinan2004}:
 
\begin{equation}\label{approx}
    \mathbb{P}(\rho_T(X^{\epsilon},\psi)<\delta) \approx \exp \left(-\frac{1}{\epsilon} S_T(\psi)\right), \quad \textrm{as } \epsilon \to 0.
\end{equation} 

There are many interpretations of the rate function in the context of statistical physics. The rate function is regarded as the free energy in macro-canonical ensemble and relative entropy in micro-canonical ensemble \cite{Touchette2009}. For more applications of the large deviation principle to statistical physics, see \cite{Touchette2009}. 

\subsection{The most likely transition path}
Let $\psi$ be the path connecting two points $x_1$ and $x_2$. As the tube size $\delta$ and parameter $\epsilon$ tends to $0$, we can regard the probability in \eqref{approx} as the transition probability of system \eqref{maineq} from the point $x_1$ to the point $x_2$ though path $\psi$ . Let $\psi^*$ be the solution of the following constrained minimization problem:
\begin{equation}\label{miniT}
    S_T(\psi^*)=\inf_{\varphi\in C_{x_1}^{x_2}([0,T])} S_T(\varphi),
\end{equation}
where $C_{x_1}^{x_2}([0,T])$ is the set of all continuous functions that start at $x_1$ and end at $x_2$. Then, as $\epsilon \to 0$, $\psi^*$ is the most likely transition path from the point $x_1$ to the point $x_2$ of system \eqref{maineq} within time $T$, in the sense that the probability that $X^\epsilon$ moves around any other path is exponentially smaller in $\epsilon$ \cite{weinan2004}. Moreover, if $S_T$ achieves its minimum value on path $\varphi^*$ for all $T>0$, that is  
\begin{equation}\label{mini}
    S_{T^*}(\varphi^*)=\inf_{T>0} \inf_{\varphi\in C_{x_1}^{x_2}([0,T])}S_T(\varphi), \quad T^* \in (0,\infty],
\end{equation}
we regard $\varphi^*$ as the most likely transition path from the point $x_1$ to  the point $x_2$. Both $\psi^*$ and $\varphi^*$ are called the instanton in physical literature related to the path integral.

Note that, the big difference between Gaussian and non-Gaussian cases is the rate function usually has an explicit form with respect to the paths for the Gaussian case, while there is no such explicit form for the non-Gaussian case. We will see the details in the following two remarks.

\begin{remark}[{\bf Gaussian case} ]
Consider the stochastic dynamical systems under Gaussian noise in $\mathbb{R}^d$, that is 
\begin{equation}\label{sdeBM}
    dX^\epsilon_t=b(X_t^\epsilon)dt+ \sqrt{\epsilon} \sigma(X_t^\epsilon) dB_t,
\end{equation}
where $B_t$ is a standard $d$-dimensional Brownian motion.
Suppose the drift term $b$ and the diffusion term $\sigma$ satisfy all the necessary conditions, for example \cite[Theorem 5.6.7]{DZ98}.  The associated good rate function is 

\begin{equation}\label{rtBM}
    S_T^{BM}(\varphi)=\inf_{g}\left\{\frac{1}{2}\int_0^T |g|^2 ds: \varphi(s)=\int_0^s b(\varphi) d\tau + \int_0^s \sigma(\varphi)g d\tau \textrm{ for all } s \in [0,T]\right\}.
\end{equation}
If the the matrix $a=\sigma \sigma^\top$ is uniformly positive definite, the rate function $S_T^{BM}$ has an explicit form\cite{DZ98} expressed by paths $\varphi$ : 
\begin{equation*}
    S_T^{BM}(\varphi)=\left \{\begin{array}{cl}
        \frac{1}{2}\int_0^T (\dot{\varphi}-b(\varphi))^\top a^{-1}(\dot{\varphi}-b(\varphi))ds & \textrm{ if the integral is finite,}\\
        \infty & \textrm{ otherwise.}
    \end{array} \right. 
\end{equation*}

The constrained minimization problem\eqref{miniT} for the most likely transition path $\tilde{\varphi}$ from the point $x_1$ to the point $x_2$ within time $T$ is then
\begin{equation}\label{miniBM}
    S_T^{BM}(\tilde{\varphi})=\inf_{\varphi \in C_{x_1}^{x_2}([0,T])} S_T^{BM}(\varphi).
\end{equation}
\end{remark}

\begin{remark}[{\bf Non-Gaussian case}]
For the stochastic dynamical systems under non-Gaussian L\'evy noise \eqref{maineq}, unlike the Gaussian noise case, the rate function $S_T$ in \eqref{rf} does not have such kind of explicit form.  Because we can only derive the ``average" value:
\begin{equation}
    \tilde{G}(s) \triangleq  \int_{\mathbb{R}^{d}\backslash \{0\}} z(g(s, z)-1) \nu(d z) =\dot{\varphi}(s)-b(\varphi(s)).
\end{equation}
For the same function $\tilde{G}$, there are actually multiple nonzero functions $g$ that satisfy the above equation. 
\end{remark}

In this paper, we aim to solve problem \eqref{mini} to obtain the most likely transition path connecting two points for stochastic dynamical system \eqref{maineq}. For this purpose, we need to choose a large enough time $T$ and then solve problem \eqref{miniT} for such $T$ to obtain the most likely transition path $\varphi^*$.

In practice, we need to solve the following constrained optimization problem for a given time $T$:
\begin{equation}\label{eq_g}
    \inf_{\varphi \in C_{x_1}^{x_2}([0,T])} \inf_{\varphi=F(g)} \int_{0}^{T} \int_{\mathbb{R}^{d}}(g(s, z) \ln g(s, z)-g(s, z)+1) \nu(d z) d s,
\end{equation}
where function $g$ is a positive measurable function. Actually, the double infimum in \eqref{eq_g} can be treated as the infimum over all function $g$ that has a corresponding function $\varphi$ in  $C_{x_1}^{x_2}([0,T])$, which is 
\begin{equation}\label{inf}
    \inf_{\{g: \varphi=F(g), ~ \varphi \in C_{x_1}^{x_2}([0,T])\}} \int_{0}^{T} \int_{\mathbb{R}^{d}}(g(s, z) \ln g(s, z)-g(s, z)+1) \nu(d z) d s.
\end{equation}

In the next section, we will construct a neural network to obtain the function $g$ and $\varphi$ which solve the minimization problem \eqref{inf}. 

\section{Method}\label{ML}
In this section, we will reformulate the constrained minimization problem \eqref{inf} as an optimal control problem. Then, we use two neural networks to approximate the corresponding state and control. With a proper loss function, we can obtain a most likely transition path by training such neural networks.

\subsection{Viewing as an optimal control problem}
 
The constrained minimization problem \eqref{inf} can be reformulated as the following optimal control problem: 
\begin{equation}\label{LoptC}
    \left\{\begin{array}{cl}
      \underset{g \in \mathcal{U}}{\operatorname{inf}}   & \bm{\mathcal{J}}[\varphi ;g]=\int_{0}^{T} \bm{\mathcal{L}}(g(s, \cdot ))d t+\eta(\varphi(T)),  \\
        \text { subject to } & \dot{\varphi}(t)=b(\varphi(t))+\bm{\mathcal{Q}}(g(t,\cdot)),\\
        &  \varphi(0)=x_1.
    \end{array}   \right.
\end{equation}
Here, $\varphi(t): [0,T] \to \mathbb{R}^d$ is the state and $g(t,\cdot): [0,T] \to \mathcal{B}^+(\mathbb{R}^d)$ is the control. $\mathcal{B}^+(\mathbb{R}^d)$ is the space of all positive Borel measurable functions and $\mathcal{U}$ is the space of positive measurable functions $g$'s that satisfy condition \eqref{condg}. The cost functional $\bm{\mathcal{J}}[\varphi; g]$ is composed of the running cost $\bm{\mathcal{L}}(g(s, \cdot))$ and the terminal cost $\eta(\varphi(T))$. The functional $\bm{\mathcal{L}}(g(s,\cdot)): \mathcal{B}^+(\mathbb{R}^d) \to \mathbb{R}^+$ and $\bm{\mathcal{Q}}(g(t,\cdot)): \mathcal{B}^+(\mathbb{R}^d) \to \mathbb{R}$ is defined by 
\begin{equation}\label{intL}
    \bm{\mathcal{L}}(g(s,\cdot)) = \int_{\mathbb{R}^d \backslash \{0\}} \big(g(s,z) \ln g(s,z) - g(s,z) +1 \big) \nu(dz),
\end{equation}
and 
\begin{equation}\label{intQ}
    \bm{\mathcal{Q}}(g(t,\cdot)) = \int_{\mathbb{R}^d \backslash \{0\}}(g(t,z)-1)z \nu(dz).
\end{equation}
Function $\eta$ is the terminal cost, and it is defined by 
\begin{equation}\label{Ecost}
    \begin{cases}
    \eta(x) = 0, & x=x_2, \\
    \eta(x) = \infty, & \textrm{otherwise.}
    \end{cases}
\end{equation}
The terminal cost function $\eta$ will guarantee that state $\varphi$ will start at point $x_1$ and end at point $x_2$, which agrees with the choose of path in the constrained minimization problem \eqref{inf}.The the optimal state $\varphi^*$ of problem \eqref{LoptC} is the most likely transition path from the point $x_1$ to the point $x_2$ within a given time $T$ for system \eqref{maineq}.
Unlike a common optimal control problem, the control in the problem \eqref{LoptC} is a map taking value in a function space. And for this reason, the known methods in the context of optimal control can not be directly applied.

However, for the stochastic dynamical system under Gaussian noise \eqref{sdeBM}, the constrained minimization problem \eqref{miniBM} can be reformulated as an optimal control problem as \cite[Equation (2.1)]{WKang2021} as follows:
\begin{equation}\label{BoptC}
    \begin{cases}
    \underset{g \in L^2([0,T]) }{\operatorname{inf}} & \mathcal{J}[\varphi; g]=\frac{1}{2}\int_{0}^{T} \mathcal{L}(g(t)) d t+\eta(\varphi(T)), \\ 
    \text { subject to } & \dot{\varphi}(t)=b(\varphi(t))+\sigma(\varphi(t))g(t), \\
    &  \varphi(0)=x_1.
    \end{cases}
\end{equation}
Here $\varphi(t): [0,T] \to \mathbb{R}^d$ is the state and $g(t) \in L^2([0,T])$ is the control. The cost functional $\mathcal{J}[\varphi;g]$ is composed of the running cost $\mathcal{L}(g(t))=|g(t)|^2$ and the terminal cost $\eta(\varphi(T))$. The terminal cost function $\eta$ is defined the same as \eqref{Ecost}. The optimal state $\varphi^*$ of problem \eqref{LoptC} is the most likely transition path from the point $x_1$ to the point $x_2$ within a given time $T$ for system \eqref{sdeBM}.
The problem \eqref{BoptC} can be solved directly by methods in the field of optimal control problem, for example\cite{WKang2021}.

As far as we know, the existing methods in the context of optimal control problem can not be directly applied to problem \eqref{LoptC}. Then, in the next subsection, we will introduce the neural network method for solving problem \eqref{LoptC}.
\subsection{Neural network method}

In this paper, we use the fully connected neural network to approximate functions $\varphi$ and $g$ in problem \eqref{LoptC} and \eqref{BoptC}. Then, we train the neural networks with a proper loss function to obtain the optimal state in problem \eqref{LoptC} and \eqref{BoptC}, which is the desired most likely transition path within time $T$.

Let $\mathcal{N}^{L}: \mathbb{R}^{n_0} \rightarrow \mathbb{R}^{n_L}$ be a fully connected neural network with  $ L$ layers and $n_{j}$ neurons in the $j$-th layer. 
The neural network is defined by
\begin{equation} \label{NNdef} 
\begin{aligned}
\mathcal{N}^{j}(\mathbf{x}) & =\left\{
        \begin{array}{ll}
          w^{j} \sigma\left(\mathcal{N}^{j-1}(\mathbf{x})\right)+b^{j}, & \text{ for } 2 \leq j \leq L,   \\
          w^{1} \mathbf{x}+b^{1},&  \text{ for } j = 1.   \\
             \end{array}
\right.
\end{aligned}
\end{equation}
where $w^{j}$ is the weight matrix in the $j$-th layer , $b^{j} \in \mathbb{R}^{n_{j}}$ is the bias vector in $j$-th layer, and $\sigma(\cdot)$ is the activation function which can be chosen as tanh, sigmoid and ReLU function. 

For the non-Gaussian case \eqref{LoptC}, we construct two fully connected neural networks to approximate the function $\varphi$ with input $t$, and $g$ with input $(t,z)$ in \eqref{LoptC}. We denote them as $\varphi_{NN}(t;w_{\varphi},b_{\varphi})$ and $g_{NN}(t,z;w_g,b_g)$ respectively, where $w_{\varphi}$ and $b_{\varphi}$ ($w_{g}$ and $b_{g}$ ) are the corresponding weights and biases.

On one hand, the neural networks $\varphi_{NN}(t;w_{\varphi},b_{\varphi})$ and $g_{NN}(x,t;w_g,b_g)$ should satisfy the 
constraint in problem \eqref{LoptC}, which means we need keep the following loss small enough:
\begin{equation*}
    \textrm{loss}_{\varphi} = \frac{1}{N_T}\sum_{i=1}^{N_T}\left(\dot{\varphi}_{NN}(t_i) -b(\varphi_{NN}(t_i)) - \mathbf{Int}_z\big(z(g_{NN}(t_i,z)-1)e^{-|z|^{\gamma}}\big)\right)^2 + \tau_1(\varphi_{NN}(0)-x_1)^2,
\end{equation*}
where $\tau_1$ is a sufficiently large positive constant so that the start-point constraint is satisfied. The residual points of the training data are $\{t_i\}_{i=1}^{N_T}$ which can be uniformly chosen on the interval $[0,T]$. The operator $\mathbf{Int}_z(f(z))$ denotes the numerical integration of function $f(z)$ with respect to variable $z$. In practice, we can use the complex trapezoid formula or the Monte Carlo method as the numerical integration $\mathbf{Int}_z$.

On the other hand, the neural network $g_{NN}(t,z)$ should be optimized to attain the minimum value in problem \eqref{LoptC}, that is to find the minimum value of the following loss:
\begin{equation*}
    \textrm{loss}_{g} =\frac{1}{N_T }\sum_{i=1}^{N_T}\mathbf{Int}_z\left(\big(g_{NN}(t_i,z)\ln g_{NN}(t_i,z)-g_{NN}(t_i,z)+1 \big)e^{-|z|^{\gamma}}\right)+\tau_2(\varphi_{NN}(T)-x_2)^2,
\end{equation*}
where $\tau_2$ is a sufficiently large positive constant to satisfy the end-point constraint.
In practice, we choose a bounded symmetric rectangular in $\mathbb{R}^d$ as the integration domain for \eqref{intL} and \eqref{intQ}, so that each integration outside this domain can be neglected. 

Combining these two constraints, we construct the total loss:
\begin{equation}\label{loss_L}
    \textrm{loss} =\tau \textrm{loss}_{\varphi}+ \textrm{loss}_{g},
\end{equation}
where $\tau$ is the weight to balance $\textrm{loss}_{\varphi}$ and $\textrm{loss}_g$.

For the Gaussian case \eqref{BoptC}, we define a similar loss function as follows:

\begin{equation}\label{loss_BM}
        \textrm{loss}^{BM}  = \tau^{BM} \textrm{loss}_{\varphi}^{BM}+ \textrm{loss}_{g}^{BM},
\end{equation}
where
\begin{align*}
        \textrm{loss}^{BM}_{\varphi} &= \frac{1}{N_T}\sum_{i=1}^{N_T}\big(\dot{\varphi}_{NN}(t_i) -b(\varphi_{NN}(t_i)) - \sigma(\varphi_{NN}(t_i))g_{NN}(t_i)\big)^2+\tau_1^{BM}(\varphi_{NN}(0)-x_1)^2, \nonumber\\
        \textrm{loss}^{BM}_g & = \frac{1}{N_T}\sum_{i=1}^{N_T}\big|g_{NN}(t_i)\big|^2+\tau_2^{BM}(\varphi_{NN}(T)-x_2)^2,
\end{align*}

where $\varphi_{NN}$ and $g_{NN}$ are neural networks for functions $\varphi$ and $g$ in problem \eqref{BoptC}. Constants $\tau_1^{BM}$ and $\tau_2^{BM}$ are sufficient large to ensure the start-point and end-point constraints. The positive constant $\tau^{BM}$ is the weight to balance $\textrm{loss}^{BM}_{\varphi}$ and $\textrm{loss}^{BM}_g$.

We will train neural networks $g_{NN}$ and $\varphi_{NN}$, so that the loss function \eqref{loss_L}  and \eqref{loss_BM} attains its minimum. In the next section, we will illustrate our method through two examples.

\section{Numerical Experiments}\label{NE}

In this section, we apply our method to compute the most likely transition path for the Maier Stein system under Gaussian and non-Gaussian L\'evy noise. Neural networks with 2 hidden layers and 20 neurons for each are implemented for the state $\varphi$. Neural networks that have 3 hidden layers with 128 neurons for each are used to approximate the control $g$ in the Gaussian case, and 2 hidden layers with 20 neurons for each are implemented for the control $g$ in the non-Gaussian L\'evy case. And we use an exponential function in the last layer for $g$ in the non-Gaussian L\'evy case.
For the neural networks, the weights are initialized with truncated normal distributions and the biases are initialized as zero. We use the $tanh$ activation function and the neural networks are implemented in TensorFlow. We use the automatic differentiation in TensorFlow to compute the derivative $\dot{\varphi}$. The Adam optimizer with a learning rate of $10^{-3}$ is used to solve the control problem. The training takes $200,000$ iterations for the following examples. 

\begin{example}[The Maier-Stein system under Gaussian noise \cite{Maier1993} ]\label{example1}
    Consider the Maier Stein system under Gaussian noise
    in $\mathbb{R}^2$:
    \begin{equation}\label{MSB}
        dX^\epsilon_t=b(X_t^\epsilon)dt+ \sqrt{\epsilon} dB_t,
    \end{equation}
    where
    \[
        b(x,y)=\left(\begin{array}{c}
            x-x^{3}- \beta x y^{2} \\-\left(1+x^{2}\right) y
            \end{array}\right).
    \]
    This system has two metastable points $(\pm 1,0)$ and one saddle point $(0,0)$. When $\beta=1$, the system is a gradient system with potential $V(x,y)=-\frac{1}{2}x^2+\frac{1}{4}x^4+\frac{1}{2}y^2+\frac{1}{2}x^2y^2$. We compute the most likely transition path from $(-1,0)$ to $(1,0)$ for parameters $\beta=10$ and $\beta=1$. We choose time $T=100$ and the number of residual points $N_T=10000$ with parameters $\tau^{BM}=1$ and $\tau_1^{BM}=\tau_2^{BM}=10$, which will enable us to obtain the most likely transition path to problem \eqref{mini}. 
   
    The results are shown in Figure \ref{2d_bm}. In Figure \ref{2d_bm}(a1), we plot the most likely transition path for the Maier Stein system with parameter $\beta=10$. 
    It shows that the most likely transition path calculated by our method is close to the one computed by the aMAM. We plot the value of loss functions with respect to the iteration step in Figure \ref{2d_bm}(a2). The value of the $\textrm{loss}_\varphi^{BM}$, which ensures the state satisfy the constraints, decreases to around $0.01$ and stays at this value. Meanwhile, the value of the $\textrm{loss}_g^{BM}$, which computes the cost, continues to decrease to near $10^{-11}$. 
    In Figure \ref{2d_bm}(b1), we plot the most likely transition path for the Maier Stein system with parameter $\beta=1$ and the value of loss functions with respect to the iteration steps in Figure \ref{2d_bm}(b2). The most likely transition path is a line between the two metastable points and the loss function also achieves a fairly small value, which also indicates the validation of our proposed method. 
 
\begin{figure}[htbp]
\begin{minipage}[]{0.5 \textwidth}
    \centerline{\includegraphics[width=8.9cm,height=7cm]{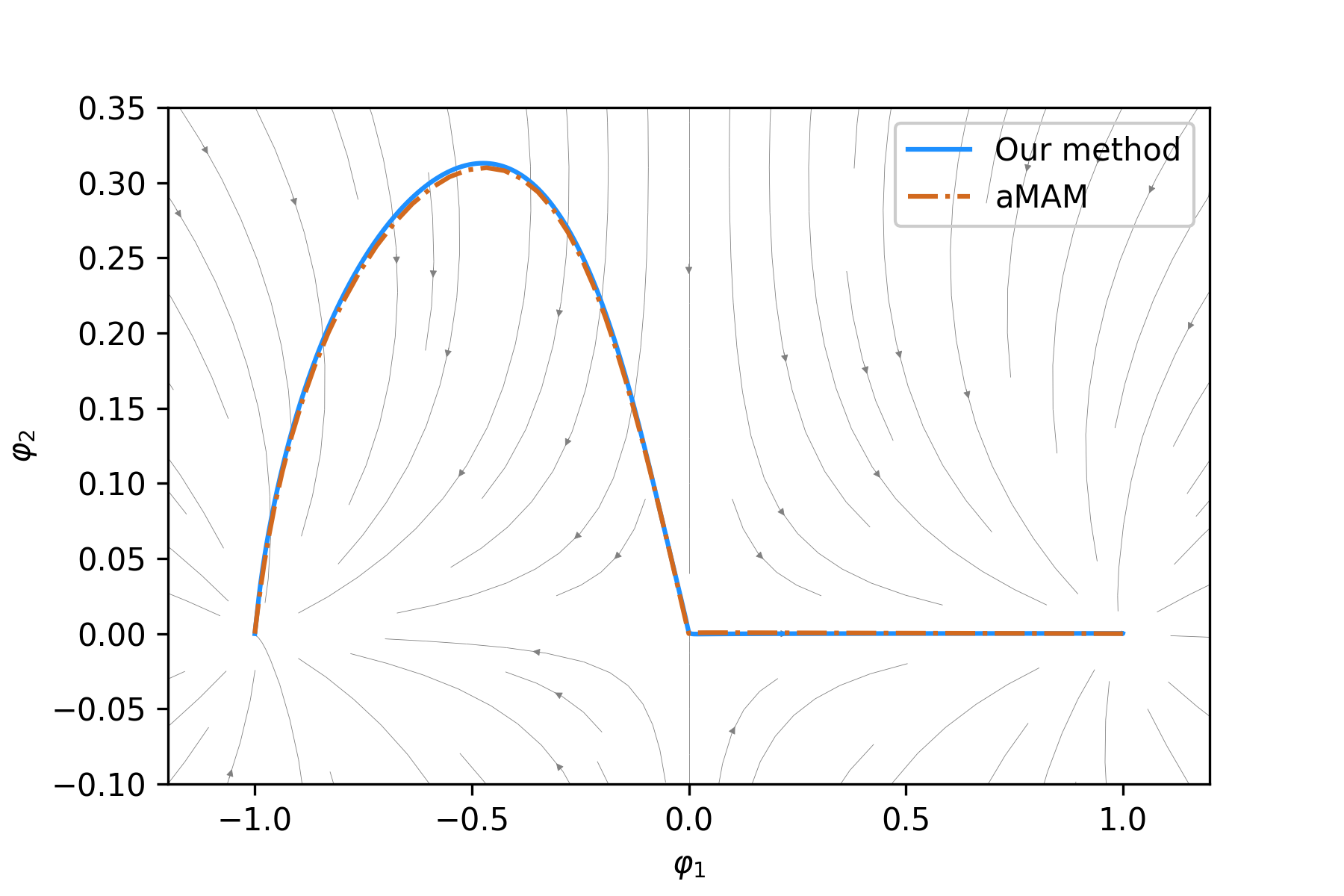}}
    \centerline{(a1)}
\end{minipage}
\hfill
\begin{minipage}[]{0.5 \textwidth}
    \centerline{\includegraphics[width=8.9cm,height=7cm]{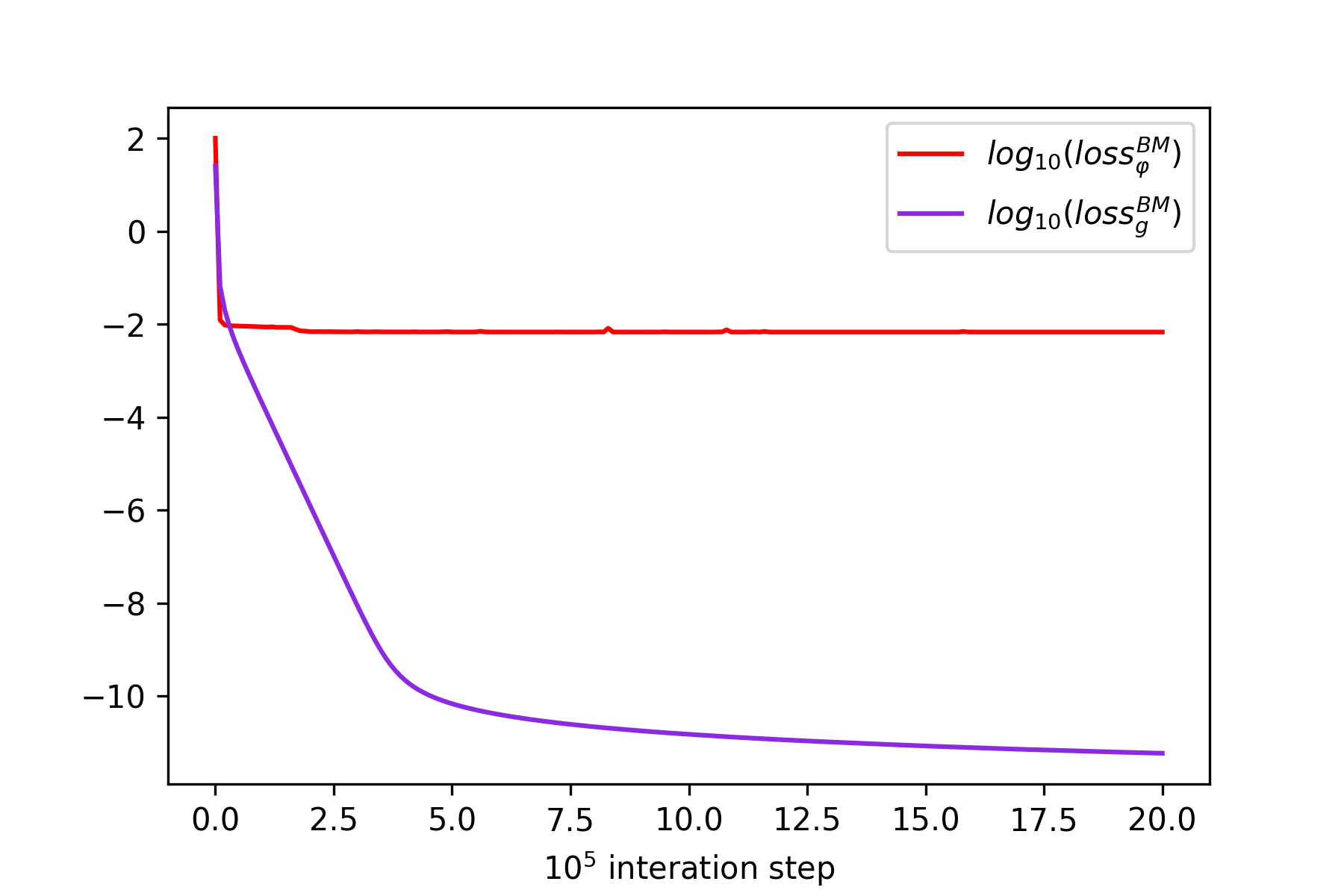}}
    \centerline{(a2)}
 \end{minipage} \\
 \begin{minipage}[]{0.5 \textwidth}
    \centerline{\includegraphics[width=8.9cm,height=7cm]{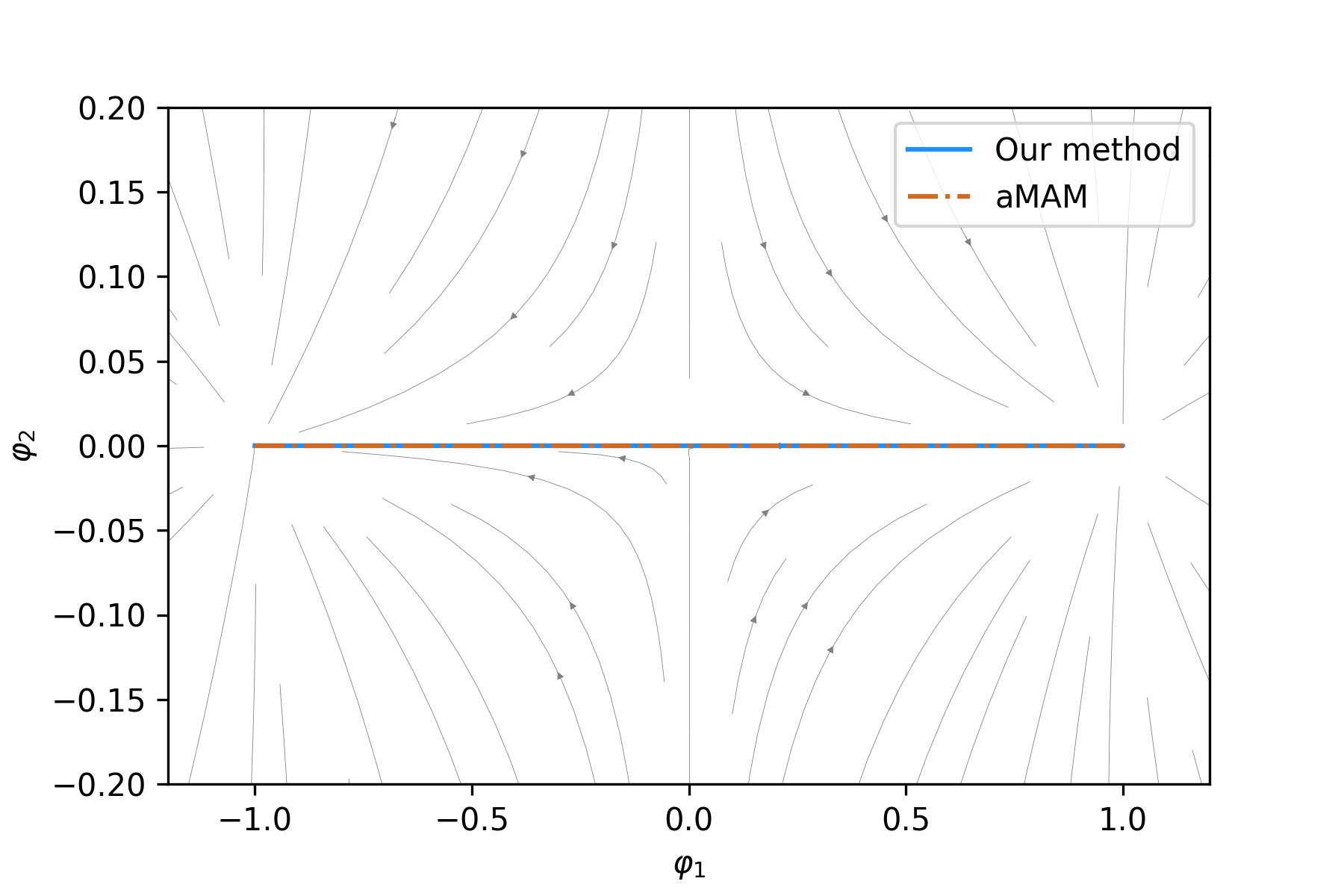}}
    \centerline{(b1)}
 \end{minipage}
 \hfill
\begin{minipage}[]{0.5 \textwidth}
    \centerline{\includegraphics[width=8.9cm,height=7cm]{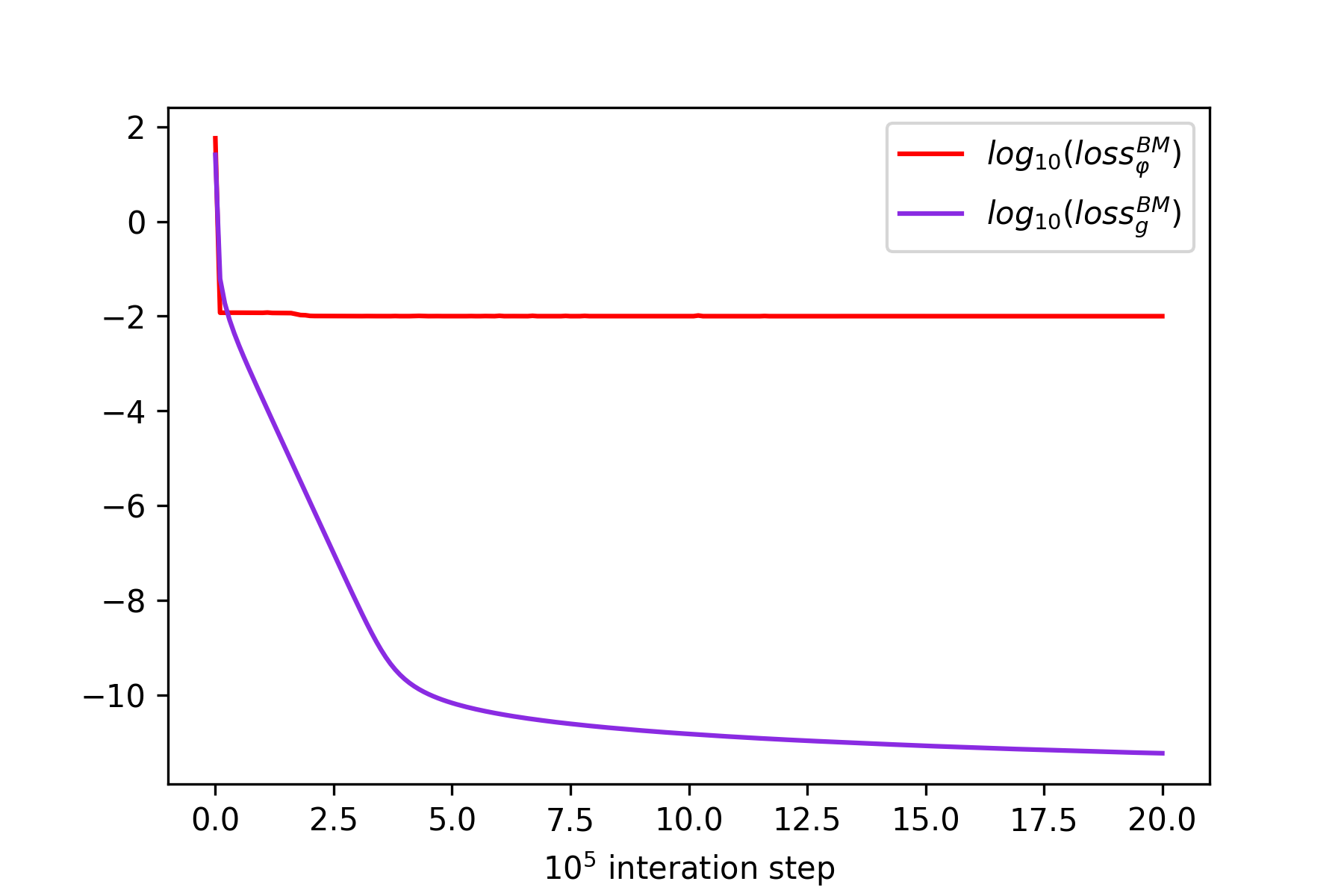}}
    \centerline{(b2)}
\end{minipage}
 
 \caption{\textbf{The case under Gaussian noise (Example \ref{example1}}): The most likely transition path from the metastable point $(-1,0)$ to the metastable point $(1,0)$  within time $T=100$ for different parameter $\beta$, computed by our method(blue line) and an adaptive minimum action method(aMAM) (brown line) on the top of the flow lines of the deterministic velocity field(grey line). (a1)-(a2) $\beta=10$: (a1) The most likely transition path; (a2) the loss function with respect to the evolution of iteration step. (b1)-(b2) $\beta=1$: (b1) the most likely transition path; (b2) the loss function with respect to the evolution of iteration step.}
 \label{2d_bm}
 \end{figure}
\end{example}

\begin{example}[The Maier-Stein system under non-Gaussian L\'evy noise]\label{example2}

\begin{figure}[htbp]
\begin{minipage}[ ]{0.5 \textwidth}
    \centerline{\includegraphics[width=8.9cm,height=7cm]{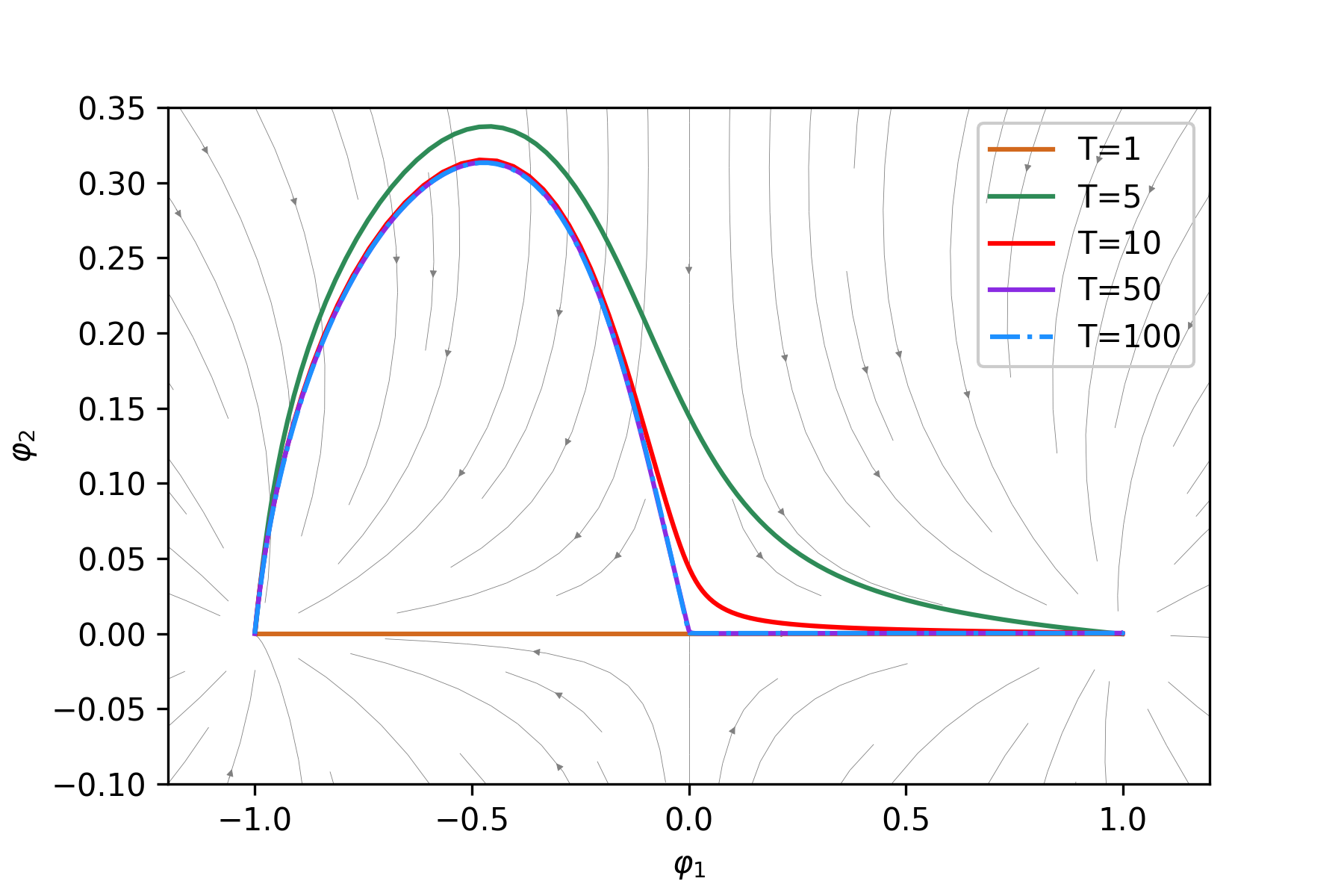}}
    \centerline{(a)}
\end{minipage}
 \hfill
\begin{minipage}[]{0.5 \textwidth}
    \centerline{\includegraphics[width=8.9cm,height=7cm]{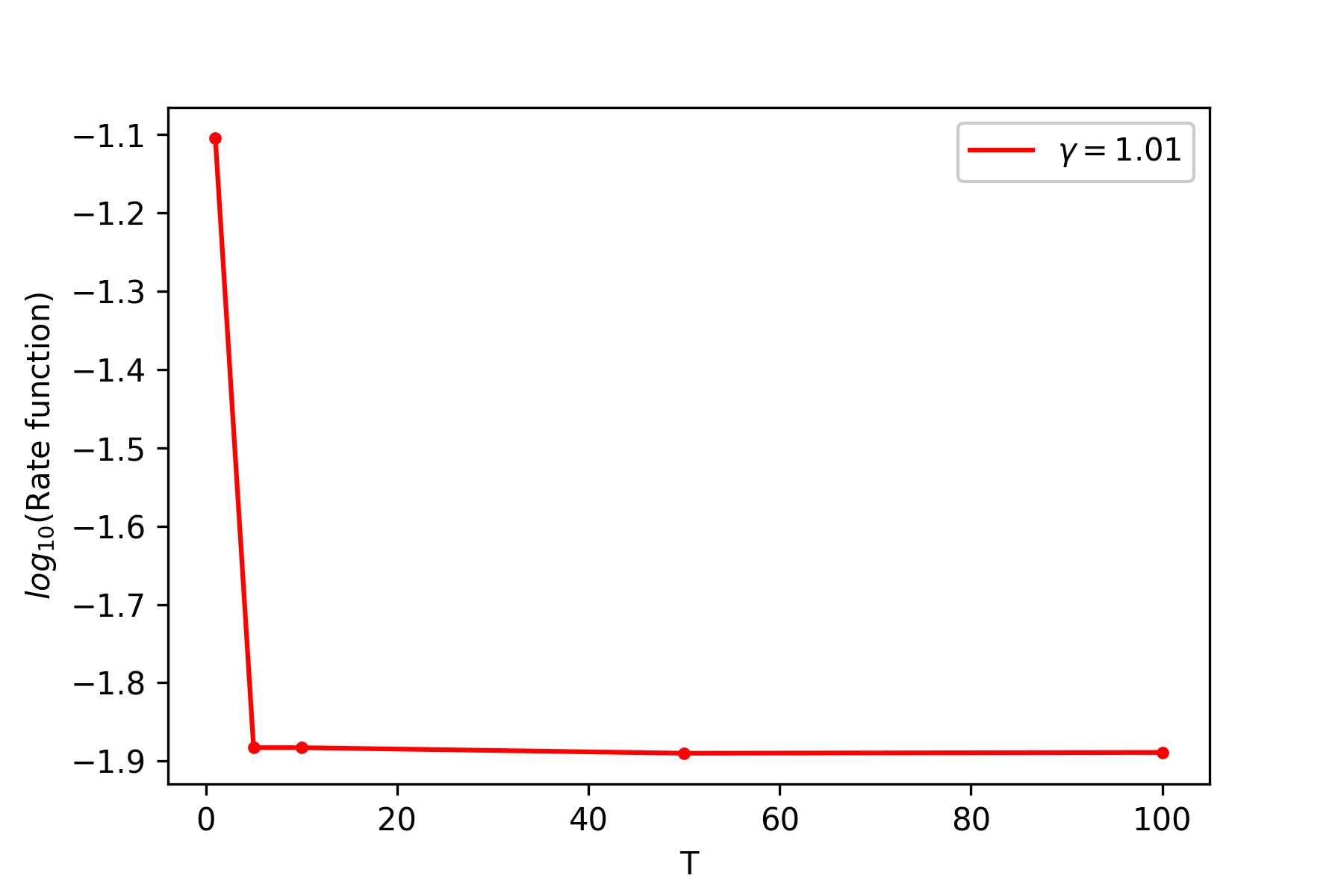}}
    \centerline{(b)}
 \end{minipage}
 \caption{\textbf{The case under non-Gaussian L\'evy noise for $\beta=10$ ( Example \ref{example2}): } (a) The most likely transition path for the Maier Stein system under non-Gaussian L\'evy noise from the metastable point $(-1,0)$ to the metastable point $(1,0)$ within time $T=1, 5, 10, 50, 100$. The grey lines are the flow lines of the deterministic velocity field. (b) The value of rate function for the most likely transition path within the time $T=1, 5, 10, 50, 100$.}\label{2d_levy}
\end{figure}
    Consider the Maier Stein system under non-Gaussian L\'evy noise \cite{YangLi2019}:
    \begin{equation}\label{MSL}
     d X_t^{\varepsilon}=b\left(X_t^{\varepsilon}\right) d t+\varepsilon d \tilde{L}_t^{\varepsilon},\quad X_0^\epsilon=x_1,
    \end{equation}
    where $b(x,y) $ is the same as Example \ref{example1} and $\nu(dx)=\exp(-|x|^{\gamma})dx$. It has the same two metastable points $(\pm 1,0)$ and one saddle point $(0,0)$.
   
    Specifically, the constraint in \eqref{LoptC} becomes
    \begin{align*}
        \dot{\varphi}_1(t)  =& x_1^{(1)}+ \big(\varphi_1(s)-\varphi_1^3(s)-\beta\varphi_1(s)\varphi_2^2(s) \big) \\
        & + \int_{\mathbb{R}^2 \backslash \{ 0\}}z_1\big( g(s,z_1,z_2)-1 \big)\exp\big( -(z_1^2+z_2^2)^{\frac{\gamma}{2}} \big) dz_1 dz_2 , \\
        \dot{\varphi}_2(t)  =& x_1^{(2)}- \big( 1+\varphi_1^2(s) \big) \varphi_2(s)  \\
        & +  \int_{\mathbb{R}^2 \backslash \{ 0\}} z_2 \big( g(s,z_1,z_2)-1 \big) \exp\big( -(z_1^2+z_2^2)^{\frac{\gamma}{2}} \big) dz_1 dz_2 .  
    \end{align*}
    The integral of the running cost to problem \eqref{LoptC} becomes 
    \begin{equation}\label{rcost_eq}
        \int_{0}^{T} \int_{\mathbb{R}^{2}\backslash \{0\}}\big[g(s, z_1,z_2) \ln g(s, z_1,z_2)-g(s, z_1,z_2)+1\big] \exp\big( -(z_1^2+z_2^2)^{\frac{\gamma}{2}} \big) dz_1 dz_2 d s.
    \end{equation}
    For the well trained neural network $g_{NN}$, we compute the rate function \eqref{rf}(integral of the running cost) though the following numerical formula:
    \begin{equation*}
    \frac{1}{N_T }\sum_{i=1}^{N_T}\mathbf{Int}_z\left(\big(g_{NN}(t_i,z_1,z_2)\ln g_{NN}(t_i,z_1,z_2)-g_{NN}(t_i,z_1,z_2)+1 \big)\exp(-(z_1^2+z_2^2)^\frac{\gamma}{2})\right).
    \end{equation*}
    
    \begin{figure}[htbp]
        \begin{minipage}[]{0.3 \textwidth}
            \centerline{\includegraphics[width=5.5cm,height=4.5cm]{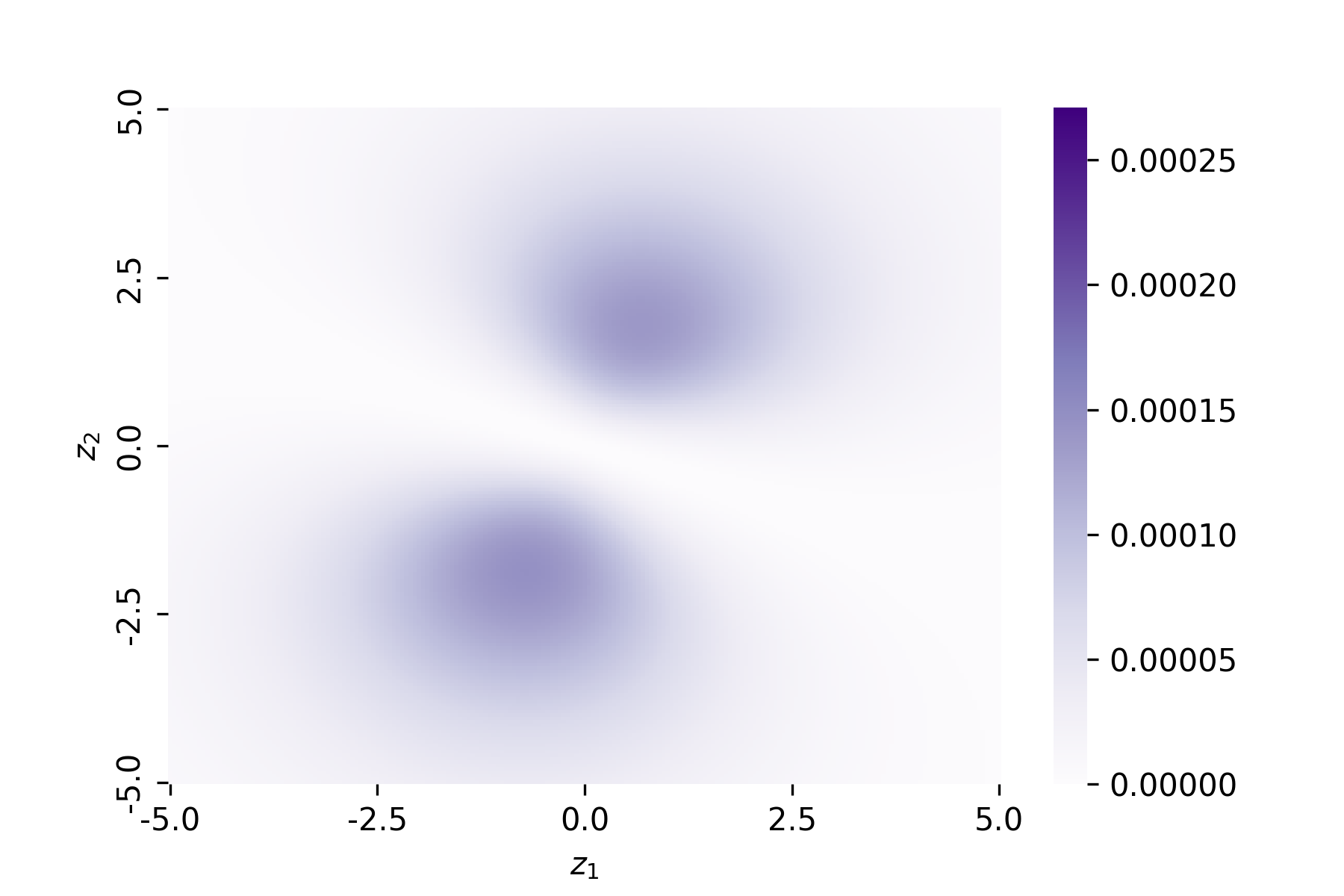}}
            \centerline{(a1)}
        \end{minipage}
        \hfill
        \begin{minipage}[]{0.3 \textwidth}
            \centerline{\includegraphics[width=5.5cm,height=4.5cm]{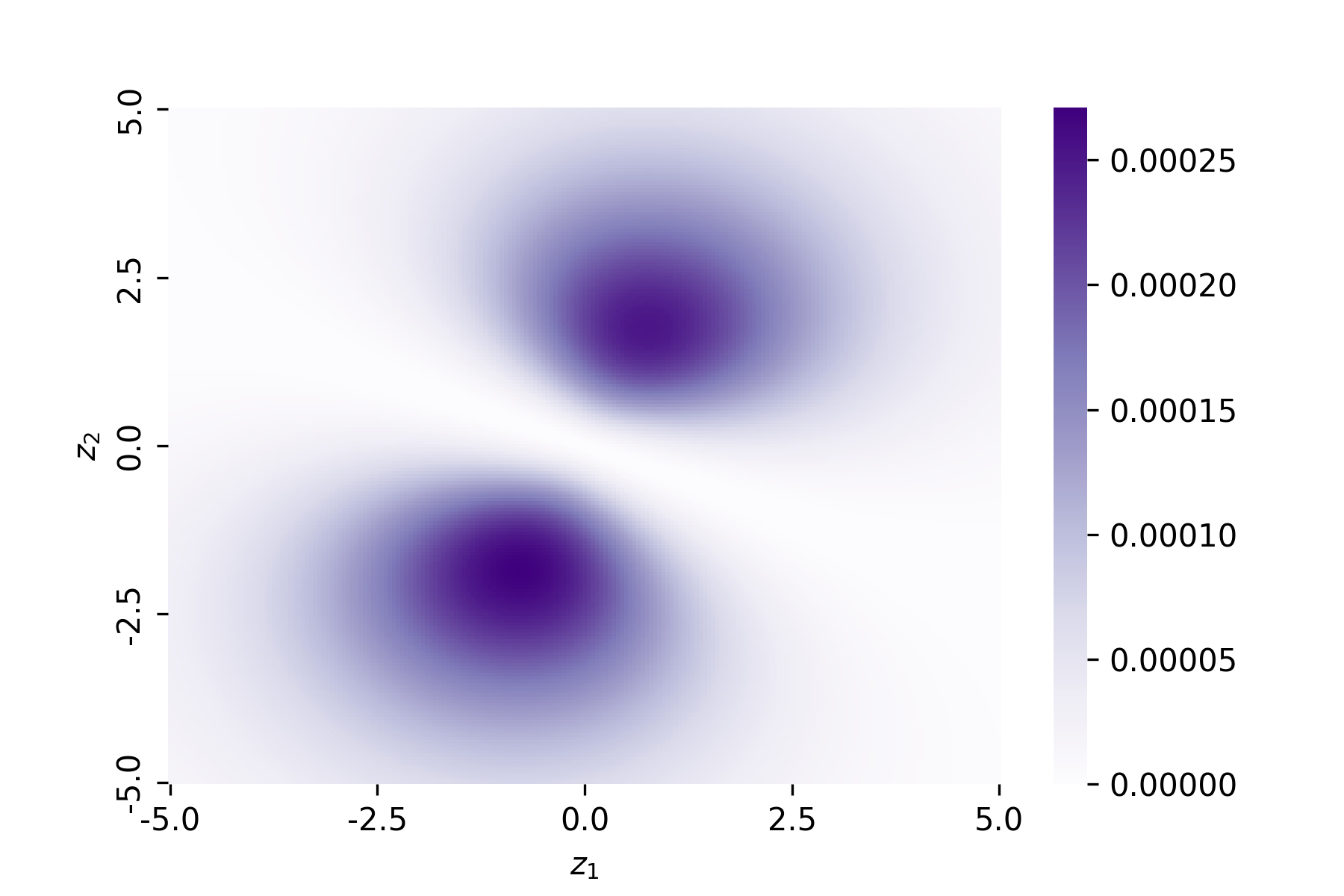}}
            \centerline{(a2)}
        \end{minipage} 
        \hfill
        \begin{minipage}[]{0.3 \textwidth}
            \centerline{\includegraphics[width=5.5cm,height=4.5cm]{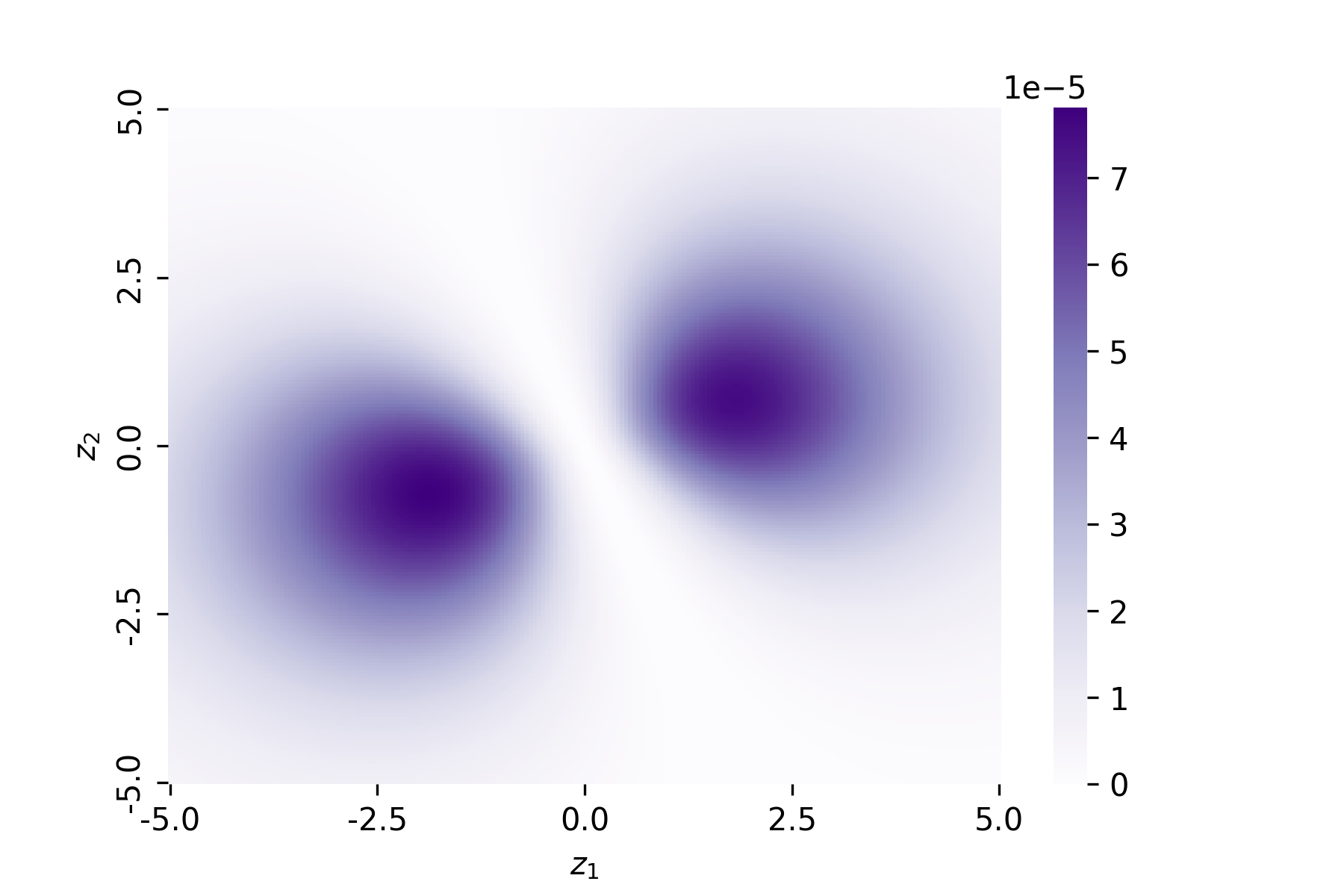}}
            \centerline{(a3)}
        \end{minipage}
        \begin{minipage}[]{0.3 \textwidth}
            \centerline{\includegraphics[width=5.5cm,height=4.5cm]{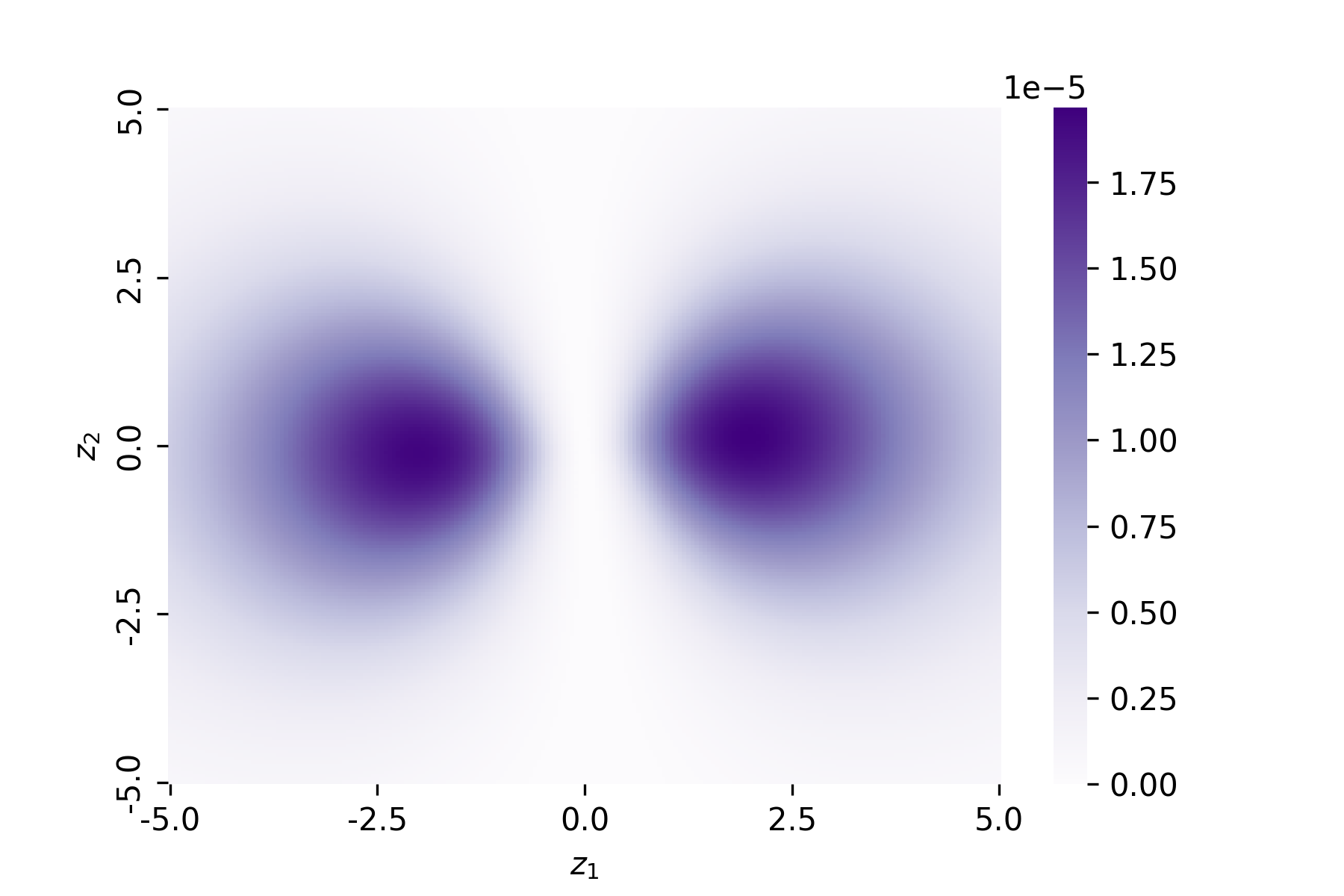}}
            \centerline{(a4)}
        \end{minipage}
        \hfill
        \begin{minipage}[]{0.3 \textwidth}
            \centerline{\includegraphics[width=5.5cm,height=4.5cm]{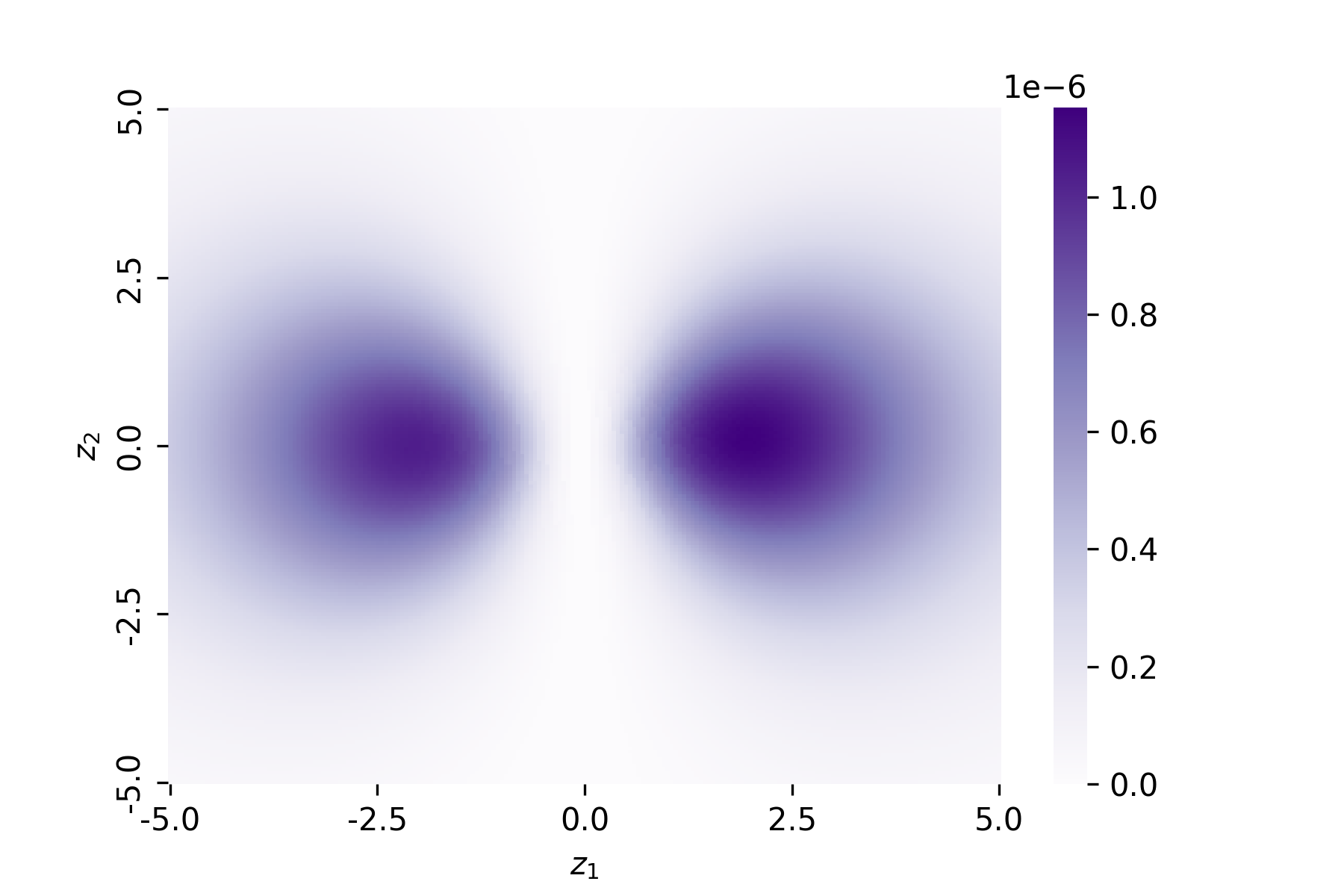}}
            \centerline{(a5)}
        \end{minipage}
        \hfill
        \begin{minipage}[]{0.3 \textwidth}
            \centerline{\includegraphics[width=5.5cm,height=4.5cm]{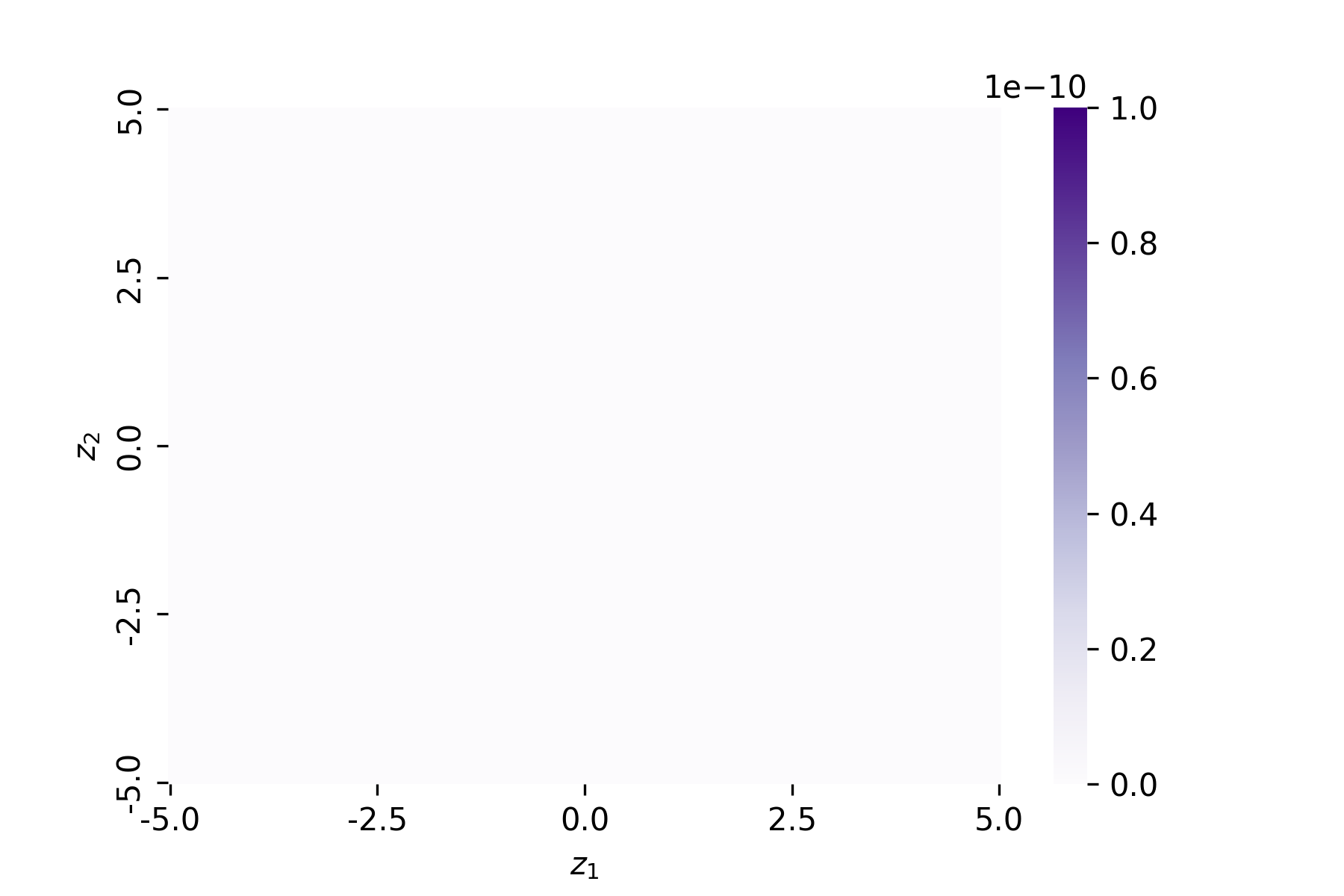}}
            \centerline{(a6)}
        \end{minipage}
        \caption{{\bf The integrand of the running cost(Example \ref{example2}):} The integrand of the running cost $(g(t,z_1,z_2)\ln g(t,z_1,z_2) -g(t,z_1,z_2) +1)\exp(-(z_1^2+z_2^2)^{\frac{\gamma}{2}})$ at different time $t$ for transition within $T=100$. (a1) $t=2$; (a2)$t=2.6$; (a3)$t=4$; (a4)$t=5$; (a5) $t=6.5$; (a6) $t=20$.}\label{costInt}
    \end{figure}
   
    We take the parameters $\beta=10$ and $\gamma=1.01$, and compute the most likely transition path from the metastable point $x_1=(-1,0)$ to the metastable point $x_2=(1,0)$ within time $T=1, 5, 10, 50$ and $100$. The number of corresponding residual points $N_T$ are chosen as $100,500,1000,1000$ and $2000$ respectively. We implement the complex trapezoid formula on the region $[-5,5]\times [-5,5]$ with uniform mesh size $\delta z_1=\delta z_2 =0.05$ as the numerical integration $\mathbf{Int_z}$, and compute with parameters $\tau=0.1$, $\tau_1=100$ and $\tau_2=10$.
    
    In Figure \ref{2d_levy}(a), the most likely transition path computed for $T=1$ is a straight line linking the two metastable points. As $T$ increases, the most likely transition path computed for $T$ converges to the one computed for $T=100$, which is similar to what is known for stochastic dynamical systems under Gaussian noise \cite{heymann2008}. The converged path(the blue dotted line in Figure \ref{2d_levy}(a)) will go through the saddle point $(0,0)$. After passing the saddle point, it follows the flow of the deterministic velocity field and reaches the metastable point $(1,0)$. 
    The value of the associated rate function is presented in Figure \ref{2d_levy}(b) with respect to time $T$. The value of the rate function also decreases to a fixed value which can be regarded as the minimum cost from $(-1,0)$ to $(1,0)$. 
    The path computed for $T=100$ is the converged path and makes the cost functional attain its minimum value. This path is then regarded as the most likely transition path from point $x_1$ to point $x_2$, which is the solution to problem \eqref{mini}.
    
    In Figure \ref{costInt}, we plot the integrand $(g(t,z_1,z_2)\ln g(t,z_1,z_2) -g(t,z_1,z_2) +1)\exp(-(z_1^2+z_2^2)^{\frac{\gamma}{2}})$ of the running cost $\bm{\mathcal{L}}$ in problem \eqref{LoptC} when $T=100$. 
    We present the integrand for time $t=2.6$ in Figure \ref{costInt}(a2), which is the one with the largest value $2.7 \times 10^{-4}$. The nonzero sets for the integrand with respect to every time $t$ stays in the region $[-5,5]\times [-5,5]$ as is shown in Figure \ref{costInt}. And it makes the numerical integration region sufficient in this case. As time $t$ increases, the value of the integrand decreases in Figure \ref{costInt}(a2)-(a6). Because the path is approaching the saddle point $(0,0)$, and less effort need the control contribute to move against the deterministic velocity field. After the path passes the saddle point, it follows a flow of the deterministic velocity field. This makes the control $g$ equal to $1$ and the integrand stays at $0$ as is shown in Figure \ref{costInt}(a6).
    
     The most likely transition paths under the non-Gaussian L\'evy case in Figure \ref{2d_levy}(a) are similar to the one under Gaussian noise in Example \ref{example1}. We also compute the most likely transition path for different $\gamma$'s. However, the computed paths are almost the same as the ones shown in Figure \ref{2d_levy}(a). The noise parameter $\gamma$ in this case does not have significant impact on the shape of the most likely transition path. It may be because the noise in the non-Gaussian L\'evy case is with exponentially light jumps and the mean square displacement \eqref{MSD} of the non-Gaussian driven noise is similar as the Gaussian driven noise. When the noise intensity $\epsilon$ tends to zero, the deterministic velocity field has the main impact on the shape of the most likely transition path as the Gaussian case, for related work, see \cite{Huang2022}.
\end{example}

\section{Discussion and Conclusion}
 In order to compute the most likely transition path for a stochastic dynamical system under non-Gaussian L\'evy noise, we have constructed a machine learning method based on an optimal control reformulation. In the context of the large deviation principle, our method finds the minimizer of the rate function among paths that connect two given points, which is regarded as the most likely transition path. Although the associated rate function in this case does not have an explicit form expressed by paths and the minimum action method can not be directly applied, our method can be still involved in.  Numerical experiments on the Maier Stein system under Gaussian and non-Gaussian L\'evy noise have demonstrated the validation of our method. Moreover, this method may also be applied to systems satisfying large deviation principle but with no explicit form for the rate function, for example, the distribution dependent system under non-Gaussian L\'evy noise and the stochastic differential equations driven by nilpotent multiplicative Gaussian noise. 
 
 As our method is implemented through neural networks, it benefits from the field of deep learning, which makes itself easy to use and capable of a wide variety of tasks.  However, it also takes the disadvantage that it is computationally expensive compared to traditional algorithms and it requires determination of a set of proper parameters. Despite these limitations, our method still shows its importance in finding the most likely transition path. 
 Further theoretical results on such optimal control problem may lead to physical informed numerical methods which may improve both the efficiency and accuracy of the algorithm.

 \section*{Acknowledgments} \noindent
 We would like to thank Dr. Yuanfei Huang and Dr. Qiao Huang for helpful discussions. This work is supported by the National Natural Science Foundation of China (NSFC) (Grant No. 12141107 and 12001067) .
 
 \section*{Data Availability}
 The data that support the findings of this study are openly available in GitHub at \\ \href{https://github.com/chenxiaolichen/Most-likely-transition-path-with-levy-noise}{https://github.com/chenxiaolichen/Most-likely-transition-path-with-levy-noise} \cite{code}.

\end{document}